\documentclass[12pt,oneside,english]{amsart}
\usepackage[T1]{fontenc}
\usepackage[latin9]{inputenc}
\usepackage{amsthm}
\usepackage{amstext}
\usepackage{graphicx}
\usepackage{setspace}
\usepackage{esint}

\makeatletter
\numberwithin{equation}{section}
\numberwithin{figure}{section}
\theoremstyle{plain}
\newtheorem{thm}{Theorem}
  \theoremstyle{definition}
  \newtheorem{defn}[thm]{Definition}
  \theoremstyle{remark}
  \newtheorem{rem}[thm]{Remark}
  \theoremstyle{plain}
  \newtheorem{cor}[thm]{Corollary}
  \theoremstyle{plain}
  \newtheorem{lem}[thm]{Lemma}

\usepackage[T1]{fontenc}
\usepackage[latin9]{inputenc}
\usepackage{verbatim}
\usepackage{wrapfig}
\usepackage{amsthm}
\usepackage{amsmath}
\usepackage{graphicx}
\usepackage{amssymb}
\usepackage[all]{xy}
\usepackage{stmaryrd}

\newcommand{\fr}{\frac}
\newcommand{\lp}{\left(}
\newcommand{\rp}{\right)}

\newcommand{\eps}{\varepsilon}
\newcommand{\trm}[1]{\textrm{#1}}
\def\hto0{\xrightarrow{h\to 0}}

\renewcommand{\d}{\partial}

\newcommand{\defeq}{\stackrel{\rm{def}}{=}}

\newcommand{\cB}{{\mathcal B}}

\newcommand{\cE}{\mathcal E}

\newcommand{\cH}{\mathcal H}

\newcommand{\cK}{\mathcal K}

\newcommand{\cN}{{\mathcal N}}

\newcommand{\cS}{\mathcal S}

\newcommand{\cV}{{\mathcal V}}

\newcommand{\ka}{\kappa}

\newcommand{\bs}{\boldsymbol}

\newcommand{\fkM}{\mathfrak M}

\newcommand{\R}{{\mathbb R}}
\newcommand{\Z}{{\mathbb Z}}

\newcommand{\supp}{\operatorname {supp}}

\newcommand{\Spec}{\operatorname{Spec}}

\newcommand{\e}{\operatorname{e}}

\newcommand{\Id}{\operatorname{Id}}

\renewcommand{\Im }{\operatorname{Im}}

\renewcommand{\i}{\operatorname{i}}

\renewcommand{\Pr}{\operatorname{Pr}}

\makeatother

\usepackage{babel}

\begin{document}
\global\long\def\trm#1{\textrm{#1}}

\title{Exponential stabilization without geometric control}

\author{Emmanuel Schenck}

\address{Mathematics Department, Northwestern University, 2033 Sheridan Road,
Evanston IL 60208, USA.}

\email{schenck@math.northwestern.edu}
\begin{abstract}
We present examples of exponential stabilization for the damped wave
equation on a compact manifold in situations where the geometric control
condition is not satisfied. This follows from a dynamical argument
involving a topological pressure on a suitable uncontrolled set. 
\end{abstract}
\maketitle

\section{Introduction}

The stabilization problem consist in studying the decay of energy
for a dissipative wave equation on a manifold \cite{BaLeRa92,Leb93,AsLe03}.
Originally aimed at controlling the vibrations of large structures,
this study has possibly numerous applications, such as for instance
sensor and actuators placements in real control problems. Furthermore,
such questions always rely on the study of non-selfadjoint operators,
the spectrum of which still lacks a good understanding. 

A simple model for the stabilization problem is given by the damped
wave equation on a compact manifold $(M,g)$ with no boundary. If
$a\in C^{\infty}(M)$, this equation reads \begin{equation}
(\d_{t}^{2}-\Delta_{g}+2a\d_{t})u=0\,.\label{eq:DWE}\end{equation}
If $a\geq0$ is non-identically zero, the energy of the waves, given
by \[
E(u,t)=\fr12(\|\nabla u\|_{L^{2}(M)}^{2}+\|\d_{t}u\|_{L^{2}(M)}^{2})\,,\]
 is decreasing in time, and satisfies $E(u,t)\stackrel{t\to\infty}{\longrightarrow}0$.
This energy decay, controlled by some specific norm of the initial
data $\bs u_{0}\defeq(u,\d_{t}u)|_{t=0}$, is at the center of the
stabilization problem. For $s>0$ one look for a function $f_{s}(t)\geq0$
such that \[
E(u,t)\leq f_{s}(t)\|\bs u_{0}\|_{\cH^{s}}^{2}\,,\quad\|\bs u_{0}\|_{\cH^{s}}^{2}=\|u(\cdot,0)\|_{H^{1+s}(M)}^{2}+\|\d_{t}u(\cdot,0)\|_{H^{s}(M)}^{2}\,.\]
A particularly interesting situation is the exponential decay, namely
when $f_{s}(t)=C\e^{-\beta_{s}t}$ for some $\beta_{s}>0$ and $C>0$
independant of $u$. The choice of the norm $\|\cdot\|_{\cH^{s}}$
is determined by the fact that the most general initial data $(u,\d_{t}u)|_{t=0}$
for the Cauchy problem arising from \eqref{eq:DWE}, belong to $\cH^{0}\defeq H^{1}(M)\times H^{0}(M)$,
see below. The case $s>0$ simply express that the energy is controlled
{}``at a greater cost of derivatives'' on the initial conditions. 

The use of global hypotheses concerning the geometry and dynamics
of the problem together with microlocal techniques has allowed a major
breakthrough concerning the stabilization problem : in \cite{RaTa75,BaLeRa92},
it has been established an exponential decay of the energy for all
initial data, under the hypothesis of geometric control (GCH). This
hypothese can be stated as follows : \smallskip{}

\begin{onehalfspace}
\noindent \begin{flushleft}
\emph{There is some $T_{0}>0$ such that every geodesic with length
$\geq T_{0}$ meets $\supp a$. }\smallskip{}

\par\end{flushleft}
\end{onehalfspace}

Until now, GCH has been an unavoidable assumption in order to get
an exponential decay \cite{BaLeRa92,Hit03}. Conversely, the failure
of GCH implies that no exponential decay can occur for waves with
the most general initial data, i.e. for $s=0$, see \cite{Leb93}.

Recently, some efforts has been made to investigate the energy decay
when GCH fails. In \cite{BuHi07}, N. Burq and M. Hitrik have studied
the damped wave equation in partially rectangular domains in the plane,
a situation that includes for instance the Bunimovich stadium. In
this latter case, assuming only $a>0$ in the two half-discs, they
were able to show a polynomial decay of the waves : for any $s>0$,
they obtained \[
E(u,t)\leq C_{s}(\log t)^{s+2}t^{-s}\|\bs u_{0}\|_{H_{0}^{1+s}\times H^{s}}^{2}.\]
Note that GCH may fail, since some bouncing ball trajectories can
stay out of $\supp a$. In an example close to what we will develop
further, Christianson \cite{Chr_H07} obtained a sub-exponential decay
for the damped wave equation on a compact Riemannian manifold, where
GCH holds outside a neighbourhood of a closed, hyperbolic orbit. For
any $s>0$, he showed that there exists some $C_{s}>0$ such that
\[
E(u,t)\leq C_{s}\e^{-C_{s}\sqrt{t}}\|\bs u_{0}\|_{\cH^{s}}^{2}.\]

\begin{figure}[h]
\begin{centering}
\includegraphics[width=8cm]{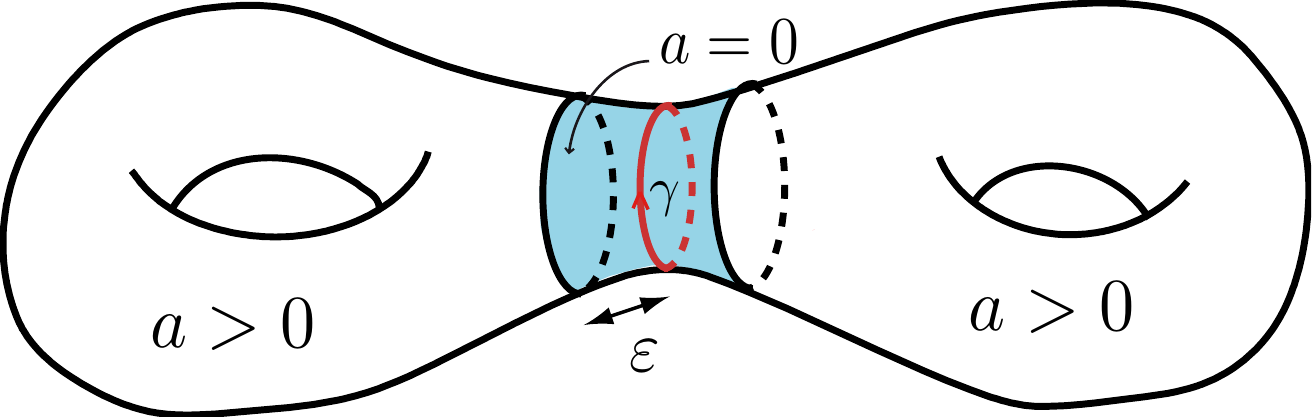}
\par\end{centering}

\caption{A simple example where the geometric control hypothesis fails : the
closed geodesic $\gamma$ is undamped. However, if $\eps$ is small
and the damping $a$ is replaced by $\beta a$ with $\beta>0$ sufficiently
large, the energy of the waves decays exponentially fast for any initial
data, in $H^{1+s}\times H^{s}$, $s>0$. \label{fig:Example}}

\end{figure}

The result we will present here address the following question : are
there situations where GCH fails, while exponential stabilization
holds ? In a recent work \cite{Sch09_2}, we showed that for manifolds
with strictly negative curvature, if GCH is dropped and replaced by
an hypothesis involving the topological pressure of the geodesic flow,
an exponential decay can happen. However, computing explicitely the
topological pressure is, in general, a difficult task, and explicit
examples showing the non-equivalence of the control and pressure hypotheses
were still lacking. As a consequence of a dynamical result of independent
interest (Theorem \ref{thm:Convergence pression} below), we provide
here a concrete example. To the author knowledge, this is the first
case where an exponential decay of the energy occurs, while the geometric
control condition fails. It can be stated as follows (see also Figure
\ref{fig:Example}) : 
\begin{thm}
\label{thm:Basic} Let $(M,g)$ be a compact Riemannian manifold of
dimension $d\geq2$, with strictly negative curvature and no boundary.
Let $\gamma\in S^{*}M$ be a periodic geodesic, and $s>0$ a real
number. For $\eps>0$, choose a function $a\in C^{\infty}(M)$, positive,
non-identically zero, with the following property : there exist a
neighbourhood $\cV$ of $\gamma$ in $S^{*}M$ and a constant $\eps$
such that

\[
\begin{array}{cl}
i) & a|_{\cV}=0\ \trm{and\ }a>0\trm{\ outside\ }\overline{\cV},\\
ii) & \trm{every\ }\rho\in\cV\ \trm{is\ at\ distance\ }\leq\eps\trm{\ from}\ \gamma\,.\end{array}\]
Then, there exists $\eps_{0}>0$ and a sequence $\{\beta_{n}\}_{n\geq0}$
satisfying $\beta_{n}\xrightarrow{n\to\infty}+\infty$ such that if
$\eps\leq\eps_{0}$ and $n$ is sufficiently large, there is $C_{s}>0$
such that for any solution $u$ of $(\d_{t}^{2}-\Delta_{g}+2\beta_{n}a\d_{t})u=0$
with initial conditions $(u,\d_{t}u)|_{t=0}=\bs u_{0}$, we have \[
E(u,t)\leq C_{s}\e^{-C_{s}t}\|\bs u_{0}\|_{\cH^{s}}^{2}\,.\]
For $s>d/2$, we have a more precise statement : for any $\epsilon>0$,
\[
E(u,t)\leq C_{\epsilon}\e^{-(G-\epsilon)t}\|\bs u_{0}\|_{\cH^{s}}^{2}\,,\]
where $G>0$ is defined in \eqref{eq:Gap}.
\end{thm}

\subsection{The topological pressure}

We will be concerned with a compact Riemannian manifold $(M,g)$ with
strictly negative curvature. We call the the geodesic flow \[
\Phi^{t}=\e^{tH_{p}}:T^{*}M\to T^{*}M\]
where \[
p:T^{*}M\owns(x,\xi)\mapsto g_{x}(\xi,\xi)=\|\xi\|_{x}^{2}\in\R\]
is the Hamiltonian and $H_{p}$ the corresponding Hamilton vector
field. Note also that $p$ is the principal symbol of $-\Delta_{g}$.
Since $M$ has strictly negative curvature, it is a well known fact
that the flow generated by $H_{p}$ on constant energy layers $\cE=p^{-1}(E)\subset T^{*}M,\ E>0$
has the Anosov property, see Appendix \ref{sec: Appendix}. The \emph{topological
pressure} with respect to the geodesic flow on $\cE$ is a functional\[
\Pr:f\in C^{\infty}(\cE)\to\overline{\R}\]
 which evaluates in some sense the complexity of the flow via the
ergodic means of the function $f$. For Anosov flows, a rather simple
definition can be given in terms of periodic orbits. A definition
suitable for general dynamical systems is given in Appendix \ref{sec: Appendix}.
For $T>0,$ let $P(T)$ be the set of closed geodesics with length
$\leq T$. Then, one has :\[
\Pr(f)\defeq\lim_{T\to\infty}\frac{1}{T}\log\lp\sum_{\gamma\in P(T)}\e^{\int_{\gamma}f}\rp\,.\]
If $f\equiv0$, one simply counts in the sum the number of periodic
orbit of length $\leq T$, and $\Pr(0)$ is equal to the topological
entropy of the flow. From now on, we restrict ourselves to the flow
on the unit cotangent bundle $\Phi^{t}:S^{*}M\to S^{*}M$, where $S^{*}M\subset T^{*}M$
is made of points with norm one covectors :\[
S^{*}M\defeq p^{-1}(1)=\{(x,\xi)\in T^{*}M:\|\xi\|_{x}=1\}\,.\]
For simplicity we will note $X\defeq S^{*}M$ in the following. This
is a compact, $2d-1$ dimensional manifold. We will consider it as
a Riemannian space too, equiped with a metric coming from $g$ such
as in the Sasaki construction. 

As shown by Walters in general settings, the topological pressure
has deep connections with invariant measures on $X$. A probability
measure $\mu$ on $X$ is said to be invariant, if for any Borel subset
$B\subset X$ and $t\geq0$ one has $\mu(\Phi^{-t}(B))=\mu(B)\,.$
We will denote by $\fkM$ the set of invariant measures on $X$. For
instance, the Liouville measure coming from the symplectic structure
of $X$ belongs to $\fkM$. The \emph{variational principle} \cite{Wal76}
asserts that one can recover the topological pressure from the knowledge
of the set $\fkM$ : \begin{equation}
\Pr(f)=\sup_{\mu\in\fkM}(h_{KS}(\mu)+\int_{X}fd\mu)\,.\label{eq:Variational Principle}\end{equation}
Here, $h_{KS}$ stands for the Kolmogorov-Sinai entropy of $\mu$.
In some cases, the $\sup$ in the above equation is attained for some
probability measure $\mu_{f}\in\fkM$ : we say that $\mu_{f}$ is
an \emph{equilibrium state} for the potential $f$. For Anosov diffeomorphisms,
such an equilibrium state always exists and is unique. 

Below, we will use the notion of topological pressure of functions
on sets $K\subset X$ which are compact and invariant under the geodesic
flow. Pesin and Pitskel \cite{PePi84} have shown that the pressure
of $f$ on $K$, denoted by $\Pr_{K}(f)$, is also given by a variational
principle. Denote by $\fkM(K)\subset\fkM$ the set of invariant measures
$\mu$ satisfying $\mu(K)=1$. Then, 

\begin{equation}
\Pr_{K}(f)=\sup_{\mu\in\fkM(K)}(h_{KS}(\mu)+\int_{K}fd\mu)\,.\label{eq:Var Pple Cpact sets}\end{equation}

\subsection{Energy decay under a pressure condition\label{sub:Energy decay}}

Before stating our main result, we recast the damped wave equation
in a standard way convenient for spectral analysis. With $\bs u=(u,\i\d_{t}u)$
one can rewrite \eqref{eq:DWE} as \begin{equation}
\d_{t}\bs u=-\i\cB\bs u,\quad\cB=\left(\begin{array}{cc}
0 & \Id\\
-\Delta_{g} & -2\i a\end{array}\right)\,.\label{eq:DWE-Matrix}\end{equation}
Via the evolution group $\e^{-\i t\cB}$, we can identify any solution
of \eqref{eq:DWE-Matrix} with its initial data $\bs u_{0}\in H^{1}\times H^{0}$.
The metastable modes of the damped wave equation \eqref{eq:DWE} can
be written in the form $\psi_{n}(t,x)=v(x)\e^{-it\tau_{n}}$ with
$\tau_{n}\in\Spec\cB$ . The spectrum of $\cB$ is a discete set of
complex numbers, localized below the real axis, in the strip $\R+\i[-2\sup a,0]$.
The only real eigenvalue is $\tau=0$, associated to constant solutions.
Let us define the function $a^{u}$ by \[
a^{u}:\begin{cases}
X\to X\\
\rho\mapsto-\bar{a}+\fr12\log J^{u}\end{cases}\]
where $J^{u}$ is the \emph{unstable Jacobian} (see Appendix \ref{sec: Appendix}),
and \[
\bar{a}:X\owns\rho\mapsto\int_{0}^{1}a\circ\Phi^{t}(\rho)dt\,.\]
The function $a$ is considered above as a function on $T^{*}M$,
depending only on the position variables. The main results of \cite{Sch09_2}
can be stated as follows :
\begin{thm}
\cite{Sch09_2} Let $(M,g)$ be a compact Riemannian manifold with
negative curvature and no boundary, and $a\in C^{\infty}(M)$ positive,
non-identically zero. Suppose that $\Pr(a^{u})<0$. Then,\begin{equation}
G\defeq\min(|\Pr(a^{u}),\inf_{\tau\in\Spec\cB\setminus\{0\}}|\Im\tau|)>0\,,\label{eq:Gap}\end{equation}
namely there is a spectral gap below the real axis. Furthermore, if
$s>d/2$, then for any $\eps>0$ we have \begin{equation}
E(u,t)\leq C_{\eps}\e^{-(G-\eps)t}\|\bs u\|_{H^{1+s}\times H^{s}}\,.\label{eq:Energy decay}\end{equation}

\end{thm}
With a standard interpolation argument, the last assertion implies
that for any $s>0$, there is some $C_{s}>0$ such that \begin{equation}
E(u,t)\leq C_{s}\e^{-C_{s}t}\|\bs u\|_{H^{1+s}\times H^{s}}\,.\label{eq:Energy decay, bis}\end{equation}

\subsection{Main result}

Since we are interested in situations where the geometric control
does not hold, we have to define an appropriate {}``undamped'' set
of $X$. It turns out that the set made of non-controlled geodesics,
given by \[
\cN\defeq\bigcap_{t\in\R}\{\rho\in X:a(\Phi^{t}(\rho))=0\}\,,\]
is too large for our purposes : we will rather define a non-controlled
subset of $X$ starting from the point of view of the invariant measures
$\fkM$ of $X$. For this, we introduce the following 
\begin{defn}
An invariant measure $\mu\in\fkM$ is \emph{minimizing} if it satisfies
$\int_{X}a\, d\mu=a_{0}$, where \[
\min_{\mu\in\fkM}\lp\int_{X}a\, d\mu\rp=a_{0}\,.\]

We call $\fkM_{0}\subset\fkM$ the set of minimizing measures on $M$.
In the interesting cases where GCH does not holds, it is clear that
$a_{0}=0$. Notice also that $\min_{\mu\in\fkM}(\int_{X}\bar{a}d\mu)=a_{0}$.
Let us now define the following set : \begin{equation}
\cK\defeq\overline{\bigcup_{\mu\in\fkM_{0}}\supp\mu}\,,\label{eq:Minimizing set}\end{equation}
which will be our {}``undamped set''. It is clear that $\cK$ is
a compact, flow-invariant subset of $X$. It has zero Liouville measure,
since this measure is ergodic with respect to the geodesic flow in
negative curvature. \end{defn}
\begin{rem}
\label{rem:K neq N} One has $\cK\subset\cN$ (see Lemma \ref{lem: a=00003D0 on K}
below), but $\cK\neq\cN$ in general. For instance, take $\gamma$
and $\beta$ two geodesics, with $\gamma$ periodic and $\beta$ homoclinic
to $\gamma$, which means that any point of $\beta$ belongs at the
same time to the stable and the unstable manifold of $\gamma$ (see
Appendix \ref{sec: Appendix}). Finally, chose $a$ positive outside
a small enough neighbourhood of $\gamma\cup\beta$, and $a=0$ near
$\gamma\cup\beta$. In this case it is easy to see that $\cN=\{\gamma\cup\beta\}$,
while $\cK=\{\gamma\}$. 
\end{rem}
Our main result, although still quite general, allows to check the
hypothesis $\Pr(a^{u})<0$ to a pressure condition on $\cK$ only,
which may be easier to evaluate : 
\begin{thm}
\label{thm:Convergence pression}Let $\bar{a}$, $a_{0}$, $J^{u}$
and $\cK$ be as above. There is a sequence $\{\beta_{n}\}_{n\geq0}$
of positive numbers with $\lim\limits _{n\to\infty}\beta_{n}=+\infty$
such that \begin{equation}
\lim_{n\to+\infty}\Pr(-\beta_{n}\bar{a}+\fr12\log J^{u})+\beta_{n}a_{0}=\Pr_{\cK}(\fr12\log J^{u})\,.\label{eq:Limite Thermodynamique}\end{equation}

\end{thm}
Basically speaking, the preceding result says that if the damping
is strong enough, the damped trajectories have a negligible contribution
to $\Pr(a^{u})$. Remain the undamped ones, for which the pressure
of $a^{u}$ reduce to the right hand side of \eqref{eq:Limite Thermodynamique}.
From this purely dynamical result, we immediately deduce the following
application : 
\begin{cor}
\label{thm: Reduction to Undamped}Suppose that $\Pr_{\cK}(\fr12\log J^{u})<0$.
Then, for some large enough $\beta>0$, we have $\Pr(-\beta a+\fr12\log J^{u})<0$.
It follows that \eqref{eq:Energy decay}, \eqref{eq:Energy decay, bis}
hold. \end{cor}
\begin{rem}
The following question, to which we are unable to answer in full generality
at this point, would have interesting consequences :
\end{rem}
\textbf{Question} : is $\cK$ a locally maximal hyperbolic subset
of $X$ ? \medskip{}

Indeed, the hypothesis $\Pr_{\cK}(\fr12\log J^{u})<0$ would be in
this case related to the size of $\cK$: a {}``filamentary'' set
$\cK$ is more likely to satisfy $\Pr_{\cK}(\fr12\log J^{u})<0$.
For instance, if $d=2$ and if $\cK$ is a locally maximal hyperbolic
set, i.e. there is a neighbourhood $U$ of $\cK$ such that \[
\cK=\bigcap_{t\in\R}\Phi^{t}(U)\,,\]
then there is a direct relationship between the Hausdorff dimension
of $\cK$ and the pressure $\Pr_{\cK}(\fr12\log J^{u})$. More precisely
\cite{PeSa01}, \[
d_{H}(\cK)<2\Longleftrightarrow\Pr_{\cK}(\fr12\log J^{u})<0\,.\]
Hence a positive answer to the above question would simply mean that
in dimension 2, checking the Hausdorff dimension of $\cK$ (or $\cN\supset\cK$)
would be sufficient to decide wether exponential stabilization can
occur or not, provided some strong enough damping is applied (note
however that it is well known that the set $\cN$ in Remark \ref{rem:K neq N}
is not locally maximal). For instance, the above remarks imply directly
Theorem \ref{thm:Basic} for $d=2$ once it is checked that $\cK=\{\gamma\}$,
which is indeed locally maximal. Finally, in higher dimension, there
is no obvious link between the Hausdorff dimension $d_{H}(\cK)$ and
the pressure $\Pr_{\cK}(\fr12\log J^{u})$, unless some assumptions
of conformality of the flow in the stable and unstable directions
are made \cite{PeSa01}.

\section{Proof of Theorem \ref{thm:Convergence pression} \label{sec:Thermo}}

We are interested in the thermodynamic limit, namely we want to study
\[
\Pr(-\beta\bar{a}+\frac{1}{2}\log J^{u})\quad\trm{in\ the\ limit\ }\beta\to+\infty\,.\]
The variational principle will be the main tool used to get Theorem
\ref{thm:Convergence pression}. We begin with two simple lemmas. 
\begin{lem}
\label{lem: a=00003D0 on K}The function $a$ satisfies \[
a|_{\cK}=a_{0}\,.\]
\end{lem}
\begin{proof}
Without loss of generality we can assume that $a_{0}=0$. Take $\rho\in\cK$,
and suppose otherwise that there is some $\rho_{0}\in\cK$ such that
$a(\rho_{0})>0$. By continuity, we can find $\epsilon\geq0$ such
that $a>0$ in some Borel subset $B(\rho_{0},\epsilon)$ of diameter
$\epsilon$ satisfying $B(\rho_{0},\epsilon)\subset\supp\mu_{0}$
for some $\mu_{0}\in\fkM_{0}$. But then,\[
\int_{B(\rho_{0},\epsilon)}ad\mu_{0}>\mu_{0}(B(\rho_{0},\epsilon))\,\inf_{B(\rho_{0},\epsilon)}a>0\,,\]
which is a contradiction since $\mu_{0}$ is minimizing. In particular,
we deduce immediately from the flow invariance of $\cK$ that $\rho\in\cN$,
showing that $\cK\subset\cN$.\end{proof}
\begin{lem}
\label{lem:Pression sur M0} We have $\fkM(\cK)=\fkM_{0}$. It follows
that for any $f\in C^{\infty}(X)$, \[
\Pr_{\cK}(f)=\sup_{\fkM_{0}}(h_{KS}(\mu)+\int_{X}fd\mu)\,.\]
\end{lem}
\begin{proof}
If $\mu\in\fkM(\cK)$, then $\mu$ is invariant and satisfies $\supp\mu\subset\cK$.
But from Lemma \ref{rem:K neq N}, we know that $\int_{X}ad\mu=a_{0}$,
which precisely means that $\mu$ is minimizing. The inclusion $\fkM_{0}\subset\fkM(\cK)$
is obvious, and the last assertion comes from \eqref{eq:Var Pple Cpact sets}.
\end{proof}
We now prove Theorem \ref{thm:Convergence pression}. By compacity,
take a sequence $(\mu_{\beta_{n}})_{n\geq0}$ of equilibrium states
for the potential $-\beta_{n}\bar{a}+\fr12\log J^{u}$ that converges
to some limit $\mu_{\infty}$ when $\beta_{n}\xrightarrow{n\to\infty}\infty$.
We first show that $\mu_{\infty}$ is minimizing. The variational
principle states that \begin{eqnarray*}
\Pr(-\beta_{n}\bar{a}+\frac{1}{2}\log J^{u}) & = & \max_{\mu\in\fkM}\lp h_{KS}(\mu)-\beta_{n}\int_{X}\bar{a}d\mu+\fr12\int_{X}\log J^{u}d\mu\rp\\
 & \leq & -\beta_{n}a_{0}+\max_{\mu\in\fkM}\lp h_{KS}(\mu)+\fr12\int_{X}\log J^{u}d\mu\rp\\
 & \le & -\beta_{n}a_{0}+\Pr(\fr12\log J^{u})\,.\end{eqnarray*}
On the other hand, we have \begin{eqnarray}
\Pr(-\beta_{n}\bar{a}+\fr12\log J^{u}) & = & \max_{\mu\in\fkM}\lp h_{KS}(\mu)-\int_{X}(\beta_{n}\bar{a}-\fr12\log J^{u})d\mu\rp\nonumber \\
 & \geq & \max_{\mu\in\fkM_{0}}\lp h_{KS}(\mu)-\int_{X}(\beta_{n}\bar{a}-\fr12\log J^{u})d\mu\rp\nonumber \\
 & \geq & -\beta_{n}a_{0}+\Pr_{\cK}(\fr12\log J^{u})\,,\label{eq: Borne sup}\end{eqnarray}
where we have used the fact that $\int_{X}\bar{a}d\mu=a_{0}$ for
$\mu\in\fkM_{0}$ and Lemma \ref{lem:Pression sur M0}. From the above
inequalities we keep in mind that\[
\Pr_{\cK}(\fr12\log J^{u})\leq\Pr(-\beta_{n}\bar{a}+\frac{1}{2}\log J^{u})+\beta_{n}a_{0}\leq\Pr(\fr12\log J^{u})\,,\]
and dividing the preceding equation by $\beta_{n}$, we obtain\[
\frac{1}{\beta_{n}}\Pr_{\cK}(\fr12\log J^{u})\leq\frac{2h_{KS}(\mu_{\beta_{n}})+\int_{X}\log J^{u}d\mu_{\beta_{n}}}{2\beta_{n}}+a_{0}-\int_{X}\bar{a}d\mu_{\beta_{n}}\leq\frac{1}{\beta_{n}}\Pr(\fr12\log J^{u})\,.\]
Letting $n\to\infty$ in the above equation shows that $\int_{X}\bar{a}d\mu_{\beta_{n}}\xrightarrow{n\to\infty}a_{0}$
, namely $\mu_{\infty}$ is minimizing. 

From \eqref{eq: Borne sup}, we have the correct lower bound in view
of Theorem \ref{thm:Convergence pression}, namely \[
\lim_{n\to\infty}\Pr(-\beta_{n}a+\fr12\log J^{u})+\beta_{n}a_{0}\geq\Pr_{\cK}(\fr12\log J^{u})\,,\]
but it remains to show the upper bound. Using the variational principle
again, we have\begin{eqnarray*}
\Pr(-\beta_{n}\bar{a}+\fr12\log J^{u}) & +\beta_{n}a_{0}= & h_{KS}(\mu_{\beta_{n}})-\beta_{n}\int_{X}(\bar{a}-a_{0})d\mu_{\beta_{n}}+\fr12\int_{X}\log J^{u}d\mu_{\beta_{n}}\\
 & \leq & h_{KS}(\mu_{\beta_{n}})+\fr12\int_{X}\log J^{u}d\mu_{\beta_{n}}\,.\end{eqnarray*}
We know that $\mu_{\beta_{n}}\xrightarrow{n\to\infty}\mu_{\infty}$,
but the upper semi-continuity of the Kolmogorov-Sinai entropy only
allows we can deduce that $\lim_{n\to\infty}h_{KS}(\mu_{\beta_{n}})\leq h_{KS}(\mu_{\infty})\,.$This
means that \begin{eqnarray*}
\lim_{\beta_{n}\to\infty}\Pr(-\beta_{n}a+\fr12\log J^{u})+\beta_{n}a_{0} & \leq & h_{KS}(\mu_{\infty})+\fr12\int_{X}\log J^{u}d\mu_{\infty}\\
 & \leq & \max_{\mu\in\fkM_{0}}(h_{KS}(\mu)+\fr12\int_{X}\log J^{u}d\mu)\\
 & \leq & \Pr_{\cK}(\fr12\log J^{u}d\mu)\,,\end{eqnarray*}
and this concludes the proof of Theorem \ref{thm:Convergence pression}.

\section{Application to the damped wave equation}

Although the Corollary \ref{thm: Reduction to Undamped} is immediate
once the limit \eqref{eq:Limite Thermodynamique} is proved, we describe
now how to obtain Theorem \ref{thm:Basic}. We use the notations of
the theorem, and show that if $\eps>0$ is small enough, then $\cK=\{\gamma\}$.
To see this, suppose that there is some $\rho\in\cV\setminus\{\gamma\}$
such that for any $t\in\R$, $\Phi^{t}(\rho)\in\cV$. Considering
a Poincaré section $S_{0}$ of the flow near $\rho_{0}\in\gamma$,
we get that under the first return map $\ka:S_{0}\to S_{0}$ the inequality
$d(\ka^{n}(\rho),\rho_{0})\leq\eps$ holds for any $n\in\Z$. But
$\ka$ being hyperbolic, it is expanding (see \cite{Bow75} for instance),
therefore we know the existence of an $\eps_{0}>0$ such that the
only points $\rho\in S_{0}$ which satisfy $d(\ka^{n}(\rho),\rho_{0})\leq\eps_{0}$
for all $n$ actually satisfy $\rho=\rho_{0}$. Hence if we take $\eps\leq\eps_{0}$,
it follows that $\cK=\{\gamma\}$. 

Let us now prove that $\Pr_{\gamma}(\fr12\log J^{u})<0$. This is
can be quickly seen from the variational principle. Observe first
that $\fkM(\gamma)=\mu_{\gamma}$ where $\mu_{\gamma}$ denotes the
invariant measure supported on $\gamma$. This implies that \begin{eqnarray*}
\Pr_{\gamma}(\fr12\log J^{u}) & = & \int_{\gamma}\fr12\log J^{u}d\mu_{\gamma}\leq-(d-1)\lambda<0\,,\end{eqnarray*}
 since $h_{KS}(\mu_{\gamma})=0$ and in the adapted metric, $J^{u}(\rho)\leq\e^{-(d-1)\lambda}$.
This concludes the proof of Theorem \ref{thm:Basic}: for $\beta$
sufficiently large, $\Pr(-\beta a+\fr12\log J^{u})$ is close enough
to $\Pr_{\gamma}(\fr12\log J^{u})<0$ to be stricly negative.

\appendix

\section{Negative curvature, Anosov flows and topological pressure\label{sec: Appendix}}

\subsection{Anosov flows. }

In this paragraph we give some further information concerning the
Anosov property of the geodesic flow in negative curvature. An Anosov
flow is a flow everywhere hyperbolic, and enjoys the following properties.
We denote as above $\cE=p^{-1}(E)\subset T^{*}M,\ E>0$ an energy
layer. For any $\rho\in\cE$, the tangent space $T_{\rho}\cE$ splits
into \emph{flow}, \emph{stable} and \emph{unstable} subspaces\[
T_{\rho}\cE=\R H_{p}\oplus E^{s}(\rho)\oplus E^{u}(\rho)\,.\]
The spaces $E^{s}(\rho)$ and $E^{u}(\rho)$ are $d-1$ dimensional,
and are preserved under the flow map: \[
\forall t\in\R,\ \ d\Phi_{\rho}^{t}(E^{s}(\rho))=E^{s}(\Phi^{t}(\rho)),\quad d\Phi_{\rho}^{t}(E^{u}(\rho))=E^{u}(\Phi^{t}(\rho)).\]
Moreover, there exist $C,\lambda>0$ such that \begin{eqnarray}
i) &  & \|d\Phi_{\rho}^{t}(v)\|\leq C\e^{-\lambda t}\|v\|,\ \trm{\ for\ all\ }v\in E^{s}(\rho),\ t\geq0\nonumber \\
ii) &  & \|d\Phi_{\rho}^{-t}(v)\|\leq C\e^{-\lambda t}\|v\|,\ \trm{\ for\ all\ }v\in E^{u}(\rho),\ t\geq0.\label{eq:Unstable}\end{eqnarray}
One can show that there exist a metric on $T^{*}M$ call the \emph{adapted
metric}, for which one can takes $C=1$ in the preceding equations.
At each point $\rho$, the spaces $E^{u}(\rho)$ are tangent to the
unstable manifold $W^{u}(\rho)$, the set of points $\rho^{u}\in\cE$
such that $d(\Phi^{t}(\rho^{u}),\Phi^{t}(\rho))\xrightarrow{t\to-\infty}0$
where $d$ is the distance induced from the adapted metric. Similarly,
$E^{s}(\rho)$ is tangent to the stable manifold $W^{s}(\rho)$, the
set of points $\rho^{s}$ such that $d(\Phi^{t}(\rho^{s}),\Phi^{t}(\rho))\xrightarrow{t\to+\infty}0$. 

The above properties allow us to define now properly the unstable
Jacobian. The adapted metric on $T^{*}M$ induces a the volum form
$\Omega_{\rho}$ on any $d$ dimensional subspace of $T(T_{\rho}^{*}M)$.
Using $\Omega_{\rho}$, we can define the unstable Jacobian at $\rho$
for time $t$. Let us define the weak-stable and weak-unstable subspaces
at $\rho$ by \[
E^{s,0}(\rho)=E^{s}(\rho)\oplus\R H_{p}\,,\quad E^{u,0}(\rho)=E^{u}(\rho)\oplus\R H_{p}.\]
We set \[
J_{t}^{u}(\rho)=\det d\Phi^{-t}|_{E^{u,0}(\Phi^{t}(\rho))}=\frac{\Omega_{\rho}(d\Phi^{-t}v_{1}\wedge\dots\wedge d\Phi^{-t}v_{d})}{\Omega_{\Phi^{t}(\rho)}(v_{1}\wedge\dots\wedge v_{d})}\,,\ \ \ \ J^{u}(\rho)\defeq J_{1}^{u}(\rho),\]
where $(v_{1},\dots,v_{d})$ can be any basis of $E^{u,0}(\Phi^{t}(\rho))$.
While we do not necessarily have $J^{u}(\rho)<1$, it is true that
$J_{t}^{u}(\rho)$ decays exponentially as $t\to+\infty$.

\subsection{Topological presure : general definition. }

Let $X$ be a metric space, and $\Phi^{t}:X\to X$ a continuous one
parameter flow. We give here the definition of the topological pressure
stated with Bowen  balls, although other equivalent definitions are
possible, for instance using open covers of $X$. For ervery $\eps>0$
and $T>0$, a set $\cS\subset X$ is $(\eps,T)-$separated if $\rho,\theta\in\cS$
implies that $d(\Phi^{t}\rho,\Phi^{t}\theta)>\eps$ for some $t\in[0,T]$,
where $d$ is the distance on $X$. For $f$ continuous on $X$, set
\[
Z(f,T,\eps)=\sup_{\cS}\left\{ \sum_{\rho\in\cS}\exp\sum_{k=0}^{T-1}f\circ\Phi^{k}(x)\right\} .\]
The topological pressure $\Pr(f)$ of the function $f$ with respect
to the flow $\Phi$ is defined by \[
\Pr(f)=\lim_{\eps\to0}\limsup_{T\to\infty}\frac{1}{T}\log Z(f,T,\eps).\]

\section*{Acknowledgments}

It is a pleasure to thank Stéphane Nonnenmacher and Philippe Thieullen
for stimulating discussions which are at the origin of this work. 

\bibliographystyle{amsalpha}
\bibliography{/Users/Admin/Documents/MATHS/LaTeX/Bibliographie/ESbiblio}

\end{document}